\input ./style/arxiv-general.cfg
\documentclass[MSNbibl,number,citesort,secthm,seceqn,dvips]{arxbj}
\makeatletter
   \@ifpackageloaded{graphicx}{}{\usepackage{graphicx}}
\makeatother

\volume{22}
\issue{4}
\pubyear{2016}
\firstpage{1963}
\lastpage{1978}
\doi{10.3150/15-BEJ716}
\docsubty{FLA}

\makeatletter
\newcommand{\rrvert}{\vert}

\newcommand{\llvert}{\vert}
\newtheorem{theorem}[thm]{Theorem}
\newtheorem{corollary}[thm]{Corollary}

\newproclaim{definition}[thm]{Definition}
\newtheorem{lemma}[thm]{Lemma}
\newremark{remark}[thm]{Remark}
\newremark{remarks}[thm]{Remarks}
\makeatother

\begin{document}
\begin{frontmatter}

\title{Conditions for a L\'evy process to stay positive near 0, in
probability}
\runtitle{L\'evy process  staying positive near 0}

\begin{aug}
\author{\inits{R.A.}\fnms{Ross A.}~\snm{Maller}\corref{}\ead
[label=e1]{Ross.Maller@anu.edu.au}}
\address{School of Finance, Actuarial Studies
and Statistics, Australian National University, Canberra, ACT,\\
Australia. \printead{e1}}
\end{aug}

%
\received{\smonth{9} \syear{2014}}

%
\begin{abstract}
A necessary and sufficient condition for a L\'evy process $X$ to stay
positive, in probability, near 0, which arises in studies of Chung-type
laws for $X$ near 0, is given in terms of the characteristics of $X$.
\end{abstract}

%
\end{frontmatter}

\section{Introduction}\label{s1}

Let $(X_t)_{t\ge0}$ be a real valued L\'evy process with canonical
triplet $(\gamma,\sigma^2,\Pi)$,
thus having characteristic function $Ee^{\mathrm{i}\theta X_t}=
e^{t\Psi
(\theta)}$, $t\ge0$,
$\theta\in\mathbb{R}$, with characteristic exponent
\begin{equation}
\label{ce}
\Psi(\theta):= \mathrm{i}\theta\gamma- \frac{1} 2
\sigma^2 \theta^2 +\int_{\mathbb{R}\setminus\{0\}}
\bigl(e^{\mathrm{i}\theta x}-1-\mathrm{i}\theta x{\mathbf{1}}_{\{
|x|\le1\}} \bigr) \Pi(\mathrm{d}x).
\end{equation}
Here, $\gamma\in\mathbb{R}$, $\sigma^2\ge0$, and $\Pi$ is a Borel measure
on $\mathbb{R}\setminus\{0\}$
such that $\int_{\mathbb{R}\setminus\{0\}}(x^2\wedge1)\Pi(\mathrm
{d}x)<\infty$.

The condition
\begin{equation}
\label{bothpos}
\mathop{\liminf}_{t\downarrow0}P(X_t\le0)\wedge
P(X_t\ge0)>0
\end{equation}
was shown by Wee \cite{wei} to imply a Chung-type law at 0 for $X$.
Attention is drawn to this in a recent paper of Aurzada, D\"oring and
Savov \cite{ads}, who give extended and refined versions of the Chung law
using a quite different approach to that of Wee.
The difference between (\ref{bothpos}) and the conditions imposed by
Aurzada et al. \cite{ads} is not at all clear, though based on some examples they suggest
that theirs are weaker than (\ref{bothpos}). Our aim in this paper is
to give necessary and sufficient conditions for $X$ to stay positive
near 0,
or to stay negative near 0, and hence to characterise~(\ref{bothpos}).

%


We need some more notation. The positive, negative and two-sided tails
of $\Pi$ are
\begin{eqnarray}
\overline{\Pi}^+(x)&:=& \Pi\bigl\{(x,\infty)\bigr\},\qquad \overline{
\Pi}^-(x):= \Pi\bigl\{(-\infty ,-x)\bigr\}\quad \mbox{and}
\nonumber
\\[-8pt]
\label{pidef}
\\[-8pt]
\nonumber
\overline{\Pi}(x) &:=&
\overline{\Pi}^+(x)+\overline{\Pi }^-(x), \qquad x>0.
\end{eqnarray}
The restriction of $\Pi$ to $(0,\infty)$ is denoted by $\Pi^{(+)}$,
and we define $\Pi^{(-)}$ on $(0,\infty)$
by $\Pi^{(-)}(\mathrm{d}x):= -\Pi(-\mathrm{d}x)$, for $x>0$.
We are only interested in small time behaviour of $X_t$, and we
eliminate the compound Poisson case by assuming $\Pi(\mathbb
{R})=\infty$ throughout.

Define truncated and Winsorised moments as
\begin{eqnarray}
\nu(x)&=& \gamma-\int_{x<|y|\le1}y\Pi(\mathrm{d}y),
\nonumber
\\[-8pt]
\label{Adef}
\\[-8pt]
\nonumber
 A(x)&=&
\gamma+\overline{\Pi}^+(1)-\overline{\Pi}^-(1) -\int_x^1
\bigl(\overline{\Pi}^+(y)-\overline{\Pi}^-(y)\bigr)\,\mathrm{d}y
\end{eqnarray}
and
\begin{equation}
\label{Vdef}
V(x)=\sigma^2+\int_{0<|y|\le x}y^2
\Pi(\mathrm{d}y),\qquad U(x)=\sigma^2+2\int_0^xy
\overline{\Pi}(y)\,\mathrm{d}y,\qquad x>0.
\end{equation}
These functions are defined and finite for all $x>0$ by virtue of
property $\int_{0<|y|\le1}y^2\Pi(\mathrm{d}y)<\infty$ of the L\'evy
measure $\Pi$ but only their behaviour as $x\downarrow0 $ will be relevant
for us.
Integration by parts shows that
\begin{equation}
\label{IP}
A(x)= \nu(x)+x\bigl(\overline{\Pi}^+(x)-\overline{\Pi}^-(x)\bigr),
\qquad x>0.
\end{equation}

Doney \cite{doney2004}, Lemma~9,  gives the following version of the It\^o
decomposition of $X$ which caters for positive and negative jumps separately.
Take constants $h_+>0$ and $ h_->0$. Then for $t\ge0$,
\begin{eqnarray}
X_t& =& t\gamma-t\nu_+(h_+) +t\nu_-(h_-)
\nonumber
\\[-8pt]
\label{decomp3}
\\[-8pt]
\nonumber
&&{}+\sigma
Z_t+ X_t^{(S,h_+,+)}+ X_t^{(S,h_-,-)}+
X_t^{(B,h_+,+)} + X_t^{(B,h_-,-)},
\end{eqnarray}
where $\gamma$ and $\sigma$ are as in (\ref{ce}), and the functions
$\nu_\pm$ are
\begin{equation}
\label{nupmdef} \nu_+(h_+):=\int_{(h_+,1]}x\Pi(\mathrm{d}x) \quad
\mbox{and}\quad \nu_-(h_-):=\int_{(h_-,1]}x\Pi^{(-)}(\mathrm{d}x).
\end{equation}
Again, only their behaviour for small values of $h_\pm$ will be
relevant. We can keep $ h_\pm\in(0,1)$.
Note that $\nu(x)=\gamma-\nu_+(x)+\nu_-(x)$. In (\ref{decomp3}),
$(X_t^{(S,h_+,+)})_{t\ge0}$ is a compensated sum of small
\textit{positive} jumps, that is,
\begin{eqnarray*}
X_t^{(S,h_+,+)}&=& \mathrm{a.s.}\lim_{\varepsilon\downarrow0}
\biggl( \sum_{0<s\le t}\Delta X_s
1_{\{\varepsilon<\Delta X_s\le h_+\}} -t\int_{\varepsilon<x\le h_+} x\Pi(\mathrm{d}x) \biggr),
\end{eqnarray*}
$(X_t^{(S,h_-,-)})_{t\ge0}$ is a compensated sum of small
\textit{negative} jumps,
that is,
\begin{eqnarray*}
&& X_t^{(S,h_-,-)}= \mathrm{a.s.}\lim_{\varepsilon\downarrow0}
\biggl( \sum_{0<s\le t}\Delta X_s
1_{\{-h_-\le\Delta X_s< -\varepsilon\}} -t\int_{-h_-\le x< -\varepsilon} x\Pi(\mathrm{d}x) \biggr),
\end{eqnarray*}
where the almost sure limits exist; and
$(X_t^{(B,h_\pm, \pm)})_{t\ge0}$ are the processes of positive and negative
big jumps, thus,
\begin{eqnarray*}
X_t^{(B,h_+,+)}&=& \sum_{0<s\le t} \Delta
X_s 1_{\{\Delta X_s>h_+\}
} \quad \mbox{and}\quad X_t^{(B,h_-,-)}= \sum
_{0<s\le t} \Delta X_s 1_{\{\Delta X_s<-h_-\}},\qquad t>0.
\end{eqnarray*}
Finally, $(Z_t)_{t\ge0}$ is a standard Brownian motion independent of
the jump processes, all of which are independent from each other.

To motivate our approach, we quote part of a result due to Doney
\cite{doney2004}. It gives an equivalence for $X$ to remain positive at small
times, with probability approaching 1, in terms of the functions
$A(x)$, $U(x)$ and the negative tail of $\Pi$.
The condition reflects the positivity of $X$ at small times in that the
function $A(x)$ remains positive for small values of $x$, and dominates
$U(x)$ and the negative tail of $\Pi$ in a certain way.

\begin{theorem}\label{th3.3}
Suppose $\Pi(\mathbb{R})=\infty$. 
\textup{(i)} Suppose also that $\overline{\Pi}^-(0+)>0$.
Then
\begin{equation}
\label{3.7}
\lim_{t\downarrow0}P(X_t> 0)=1
\end{equation}
if and only if
\begin{equation}
\label{3.6a}
\lim_{x\downarrow0} \frac{A(x)}{\sqrt{U(x)\overline{\Pi}^-(x)}}=\infty.
\end{equation}

\textup{(ii)} Suppose alternatively that
$X$ is spectrally positive, that is, $\overline{\Pi}^-(x)=0$ for
all $x>0$. Then~(\ref{3.7}) is equivalent to
\begin{equation}\label{3.8}
\sigma^{2}=0\quad \mbox{and}\quad A(x)\ge0 \qquad \mbox{for all   small }x,
\end{equation}
and this happens if and only if $X$ is a subordinator.
Furthermore, we then have $A(x)\geq0$, not only for small
$x$, but for all $x>0$.
\end{theorem}

\begin{remarks*}
(i) Other equivalences for (\ref{3.7}) are in
Theorem~1 of Doney \cite{doney2004} (and his remark following the theorem).
He assumes  a priori  that $\sigma^2=0$ but this is not necessary
as it follows from
the inequality:
\begin{equation}
\label{4.1x}
\mathop{\limsup}_{x\downarrow0}\frac{A(x)}{\sqrt{\overline{\Pi
}^-(x)}}<\infty,
\end{equation}
whenever $\overline{\Pi}^-(0+)>0$,
which is proved in Buchmann, Fan and Maller \cite{bfm}.

(ii) When $\Pi(\mathbb{R})<\infty$, $X$ is compound Poisson and its
behaviour near 0 is simply determined by the sign of the shift constant
$\gamma$. We eliminate this case throughout.
\end{remarks*}

The next section contains our main result which is essentially a
subsequential version of Theorem~\ref{th3.3}.

\section{Staying positive near 0, subsequential version}\label{s2a}
Denote the jump process of $X$ by $(\Delta X_t)_{t\ge0}$, where
$\Delta X_t = X_t-X_{t-}$, $t>0$, with \mbox{$\Delta X_0 \equiv0$}, and
define $\Delta X_t^+=\max(\Delta X_t,0)$, $\Delta X_t^-=\max(-\Delta X_t,0)$,
$(\Delta X^+)_t^{(1)}=\sup_{0<s\le t}\Delta X_s^+$, $(\Delta
X^-)_t^{(1)}=\sup_{0<s\le t}\Delta X_s^-$.

\begin{theorem}\label{th3.3s}
Assume $\Pi(\mathbb{R})=\infty$.
\begin{longlist}[(iii)]
\item[(i)]  Suppose $\overline{\Pi}^-(0+)>0$.
Then the following are equivalent:\\
there is a non-stochastic sequence $t_k\downarrow0$ such that
\begin{equation}
\label{3.7s}
P(X_{t_k}>0)\to1; 
\end{equation}
there is a non-stochastic sequence $t_k\downarrow0$ such that
\begin{eqnarray}
\label{Xdels}
\frac{X_{t_k}}{(\Delta X^-)_{t_k}^{(1)}} &\stackrel{\mathrm {P}} {\longrightarrow}& \infty\qquad
\mbox{as }  k\to\infty;
\\
\label{3.6as}
\limsup_{x\downarrow0} \frac{A(x)}{\sqrt{U(x)\overline{\Pi}^-(x)}} &=& \infty.
\end{eqnarray}

\item[(ii)]  Suppose alternatively that $X$ is spectrally positive, that is,
$\overline{\Pi}^-(x)=0$ for
all $x>0$. Then~(\ref{3.7s}) is equivalent to
$\lim_{t\downarrow0} P(X_t>0)\to1$, thus to (\ref{3.8}),
equivalently, $X_t$ is a subordinator, and $A(x)\geq0$ for all $x>0$.

\item[(iii)]  Suppose $\overline{\Pi}^-(0+)>0$. Then
$X_{t_k}/t_k\stackrel{\mathrm{P}}{\longrightarrow}\infty$ for a
non-stochastic sequence $t_k\downarrow0$
if and only if
\begin{equation}
\label{3.6as+}
\limsup_{x\downarrow0} \frac{A(x)}{1+\sqrt{U(x)\overline{\Pi}^-(x)}}=\infty.
\end{equation}
\end{longlist}
\end{theorem}

\begin{remarks*}
(i)
When $\overline{\Pi}^-(0+)>0$, $\sup_{0<s\le t} \Delta X_s^->0$ a.s.
for all $t>0$,
so the ratio in (\ref{Xdels}) is well defined.

(ii) Sato \cite{S99}, page 65,  shows that $P(X_t\le x)$ is a continuous
function of $x$ for all $t>0$ when $\Pi(\mathbb{R})=\infty$.
So $P(X_t>0)=P(X_t\ge0)$ for all $t>0$ and
$P(X_{t_k}>0)$ can be replaced by $P(X_{t_k}\ge0)$ in (\ref{3.7s})
without changing the result (and similarly in Theorem~\ref{th3.3}).

(iii) Assuming $\overline{\Pi}^+(0+)=\infty$ and $\overline{\Pi
}^-(0+)>0$, the
contrapositive of (\ref{3.7s}) shows that there is no sequence
$t_k\downarrow0$ such that $P(X_{t_k}>0)\to1$, or, equivalently,
$\liminf_{t\downarrow0}P(X_t\le0)>0$, if and only if
\begin{equation}
\label{3.6ask}
\limsup_{x\downarrow0} \frac{A(x)}{\sqrt{U(x)\overline{\Pi
}^-(x)}}<\infty.
\end{equation}
By a symmetrical argument, when $\overline{\Pi}^-(0+)=\infty$ and
$\overline{\Pi}
^+(0+)>0$, then $\liminf_{t\downarrow0}P(X_t\ge 0)>0$ if and only if
\begin{equation}
\label{3.6ask-}
\liminf_{x\downarrow0} \frac{A(x)}{\sqrt{U(x)\overline{\Pi}^+(x)}}>-\infty.
\end{equation}
Combining these gives the following.
\begin{corollary}\label{cor1}
Assume $\overline{\Pi}^+(0+)=\overline{\Pi}^-(0+)=\infty$.
Then (\ref{bothpos}) holds if and only if
\begin{equation}
\label{3.6askb}
-\infty< \liminf_{x\downarrow0} \frac{A(x)}{\sqrt{U(x)\overline{\Pi}^+(x)}}\quad
\mbox{and} \quad \limsup_{x\downarrow0} \frac{A(x)}{\sqrt{U(x)\overline{\Pi}^-(x)}}<\infty.
\end{equation}
\end{corollary}
When one of $\overline{\Pi}^+(0+)$ or $\overline{\Pi}^-(0+)$ is
infinite but the other
is zero, conditions for (\ref{bothpos}) can also be read from
Theorem~\ref{th3.3s}.

(ii)
A random walk version of Theorem~\ref{th3.3s} is in Kesten and Maller
\cite{KM}. Andrew \cite{PA}, Theorem~4,  has results related to Theorem~\ref
{th3.3s}, including the equivalence of (\ref{3.7s}) and (\ref{Xdels}).
\end{remarks*}

\section{Some inequalities for the distribution of $X$}\label{subs2.5a}
For the proof of Theorem~\ref{th3.3s}, some lemmas are needed.
The first gives a non-uniform Berry--Esseen bound for a small jump
component of $X$.
The proof is rather similar to that of Lemma~4.3 of
Bertoin, Doney and Maller \cite{bdm}, so we omit details.

\begin{lemma}\label{ron2lem1}
Fix $h_-\ge0$, $h_+\ge0$, $h_-\vee h_+>0$.
Let $(X_t^{(-h_-,h_+)})_{t\ge0}$ be the small jump martingale obtained
from $X$ as the compensated sum of jumps with magnitudes in
$(-h_-,h_+)$:
\begin{eqnarray*}
&&
X_t^{(-h_-,h_+)}= \mathrm{a.s.}\lim_{\varepsilon\downarrow0}
\biggl( \sum_{0<s\le t}\Delta X_s
1_{\{ \Delta X_s\in(-h_-,-\varepsilon)\cup(\varepsilon,h_+) \}} -t\int_{x\in(-h_-,-\varepsilon)\cup(\varepsilon,h_+)}x\Pi(\mathrm {d}x) \biggr).
\end{eqnarray*}
(Interpret\vspace*{1pt} integrals over intervals of the form
$(0,-\varepsilon)$, and $(\varepsilon,0)$, $\varepsilon>0$, as 0.)
Define absolute moments
$m_k^{(-h_-,h_+)}:=\int_{-h_-<x<h_+}|x|^k\Pi(\mathrm{d}x)$,
$k=2,3,\ldots,$
and assume $\sigma^2+m_2^{(-h_-,h_+)}>0$.
Then
we have the non-uniform bound: for any $x\in\mathbb{R}$, $t>0$,
\begin{eqnarray}\label{be1}
&& \biggl\llvert P \biggl( \frac{\sigma Z_t+X_t^{(-h_-,h_+)}}{\sqrt{t(\sigma
^2+m_2^{(-h_-,h_+)})}}\le x \biggr)-\Phi(x)\biggr\rrvert \leq
\frac{C m_3^{(-h_-,h_+)} }{\sqrt{t}( \sigma^2+m_2^{(-h_-,h_+)}
)^{3/2}(1+|x|)^{3}},
\end{eqnarray}
where $C$ is an absolute constant and $\Phi(x)$ is the standard normal c.d.f.
\end{lemma}

Next, we use Lemma~\ref{ron2lem1} to develop other useful bounds.
Define
\begin{equation}
\label{Vpmdef}
V_+(x)=\int_{0<y\le x} y^2 \Pi(\mathrm{d}y) \quad \mbox{and}\quad V_-(x)=\int_{-x\le y<0} y^2
\Pi(\mathrm{d}y), \qquad x>0.
\end{equation}
%
In the next lemma, the ``$+$'' and ``$-$'' signs are to be taken
together. When $\overline{\Pi}^+(0+)=0$ we have $V_+\equiv0$, and interpret
$(X_t^{(S,d_+,+)})_{t\ge0}$ as 0; similarly with ``$-$'' replacing
``$+$'' when $\overline{\Pi}^-(0+)=0$.

\begin{lemma}\label{dolem11}
\textup{(i)}
Suppose $d_\pm>0$, $\kappa_\pm>0$ and $K_\pm$ are constants satisfying
\begin{equation}
\label{Kpm}
K_\pm\ge4C\max \biggl(\frac{\kappa_\pm}{\Phi(-\kappa_\pm)},
\frac{1}{\Phi(-\kappa_\pm)\sqrt{1-\Phi(-\kappa_\pm)/2}} \biggr),
\end{equation}
where $C$ is the absolute constant in (\ref{be1}).
Then for each $t>0$
\begin{equation}
\label{an}
P \bigl(X_t^{(S,d_\pm,\pm)}\le K_\pm
d_\pm-\kappa_\pm\sqrt{tV_\pm (d_\pm)}
\bigr) \ge\Phi(-\kappa_\pm)/2.
\end{equation}

\textup{(ii)} 
Suppose, for each $t>0$, $d_\pm=d_\pm(t)>0$ satisfy
\begin{equation}
\label{Pilb}
t\overline{\Pi}^+(d_+)\le c_+ \quad \mbox{and}\quad t\overline{\Pi}^-(d_--)\ge
c_-
\end{equation}
for some $c_+>0$, $c_->0$.
Assume $\kappa_\pm>0$ and $K_\pm$ are constants satisfying
(\ref{Kpm}).
\begin{longlist}[000(a)]
\item[(a)]  Suppose $\overline{\Pi}^+(0+)>0$.
Then for each $t>0$ and $L\ge0$
\begin{eqnarray}
&& P \bigl(X_t\le t\gamma-t\nu_+(d_+)+t\nu_-(d_-)
+K_+d_+-Ld_- 
%
-\kappa_+
\sqrt{tV_+(d_+)} -\kappa_-\sqrt{tV_-(d_-)} \bigr)\qquad
\nonumber
\\[-8pt]
\label{Xlb2}
\\[-8pt]
\nonumber
&&\quad\ge e^{-c_+}\Phi(-\kappa_+) \Phi(-\kappa_-)P\bigl(N(c_-)\ge K_-+L
\bigr)/8,
\end{eqnarray}
where $N(c_-)$ is a Poisson rv with expectation $c_-$.

\item[(b)]  When $\overline{\Pi}^+(0+)=0$, (\ref{Xlb2}) remains true with
$\nu_+(d_+)=V_+(d_+)=d_+=c_+=0$.

\textup{(iii)}  Suppose $0\le\overline{\Pi}^-(0+)<\infty=\overline{\Pi}^+(0+)$
and, for $t>0$,
$d_+=d_+(t)>0$ is such that $t\overline{\Pi}^+(d_+(t))\le c_+$.
Suppose $\kappa_+>0$ and $K_+$ are constants satisfying (\ref{Kpm}).
Then 
\begin{equation}
\label{Xlb3}
P \bigl(X_t\le t\gamma-t\nu_+(d_+)+t\nu_-(0) +K_+d_+-
\kappa_+\sqrt{tV_+(d_+)} \bigr) \ge e^{-c_+} \Phi(-\kappa_+)/4,
\end{equation}
where $\nu_-(0)\equiv0$ when $\overline{\Pi}^-(0+)=0$.
\end{longlist}
\end{lemma}

\begin{pf}
(i) We give the proof just for the ``$+$'' signs.
Fix $t>0$ and take any constants $d_+>0$, $\kappa_+>0$ and $K_+$, with
$\kappa_+$ and $K_+$ satisfying (\ref{Kpm}).
\begin{longlist}[(a)]
\item[(a)]  Assume $V_+(d_+)>0$. Apply the bound (\ref{be1}) in Lemma~\ref
{ron2lem1} to $X_t^{(S,d_+,+)} $, which has L\'evy measure $\Pi$
restricted to $(0,d_+)$.
Noting that
$\int_{0 <y\le x}y^3\Pi(\mathrm{d}y) \le xV_+(x)$, $x>0$, (\ref{be1})
then gives, for each $t>0$,
\begin{equation}
\label{norbd}
\sup_{x\in\mathbb{R}}\bigl\llvert P \bigl(X_t^{(S,d_+,+)}
\le x\sqrt{tV_+(d_+) } \bigr) -\Phi(x)\bigr\rrvert 
\leq
\frac{Cd_+}{\sqrt{tV_+(d_+)}}.
\end{equation}
Substitute $x=-\kappa_+$ in this to get
\[
P \bigl(X_t^{(S,d_+,+)}\le-\kappa_+\sqrt{tV_+(d_+)} \bigr) \ge
\Phi(-\kappa_+) -\frac{Cd_+}{\sqrt{ tV_+(d_+)}}.
\]
When $2Cd_+\le\Phi(-\kappa_+)\sqrt{ tV_+(d_+)}$, this inequality implies
\begin{equation}
\label{ac0}
P \bigl(X_t^{(S,d_+,+)}\le-\kappa_+
\sqrt{tV_+(d_+)} \bigr) \ge\tfrac{1}{2} \Phi(-\kappa_+).
\end{equation}
When $2Cd_+> \Phi(-\kappa_+)\sqrt{ tV_+(d_+)}$,
we have
\begin{eqnarray*}
2\kappa_+\sqrt{tV_+(d_+)}< 4Cd_+\kappa_+/\Phi(-\kappa_+)\le K_+d_+,
\end{eqnarray*}
since $K_+$ satisfies (\ref{Kpm}).
Apply Chebychev's inequality, noting that
$X_t^{(S,d_+,+)}$ has mean 0 and variance $tV_+(d_+)$,
to get
\begin{eqnarray*}
P \bigl(X_t^{(S,d_+,+)}\le K_+d_+-\kappa_+\sqrt{tV_+(d_+)}
\bigr) &\ge& 1-\frac{tV_+(d_+)}{ (K_+d_+-\kappa_+\sqrt
{tV_+(d_+)} )^2}
\nonumber
\\
&\ge& 1- \frac{4tV_+(d_+)}{K_+^2d_+^2}.
\end{eqnarray*}
Also when $2Cd_+>\Phi(-\kappa_+)\sqrt{ tV_+(d_+)}$, by choice of
$K_+$ in (\ref{Kpm}) we have
\[
\frac{4tV_+(d_+)}{K_+^2d_+^2} \le\frac{16C^2}{\Phi^2(-\kappa_+)K_+^2} \le1-\frac{\Phi(-\kappa_+)}{2},
\]
giving
\begin{equation}
\label{ac2}
P \bigl(X_t^{(S,d_+,+)}\le K_+d_+-\kappa_+
\sqrt{tV_+(d_+)} \bigr) \ge\tfrac{1}{2} \Phi(-\kappa_+).
\end{equation}
The same inequality holds when $2Cd_+\le\Phi(-\kappa_+)\sqrt{
tV_+(d_+)}$, by (\ref{ac0}), so it holds in general.

\item[(b)]  When $V_+(d_+)=0$, $\Pi(\cdot)$ has no mass in $(0,d_+)$, and
(\ref{an}) with a ``$+$'' sign
remains valid in the sense that $X_t^{(S,d_+,+)}=0$ a.s.
and the left-hand side of (\ref{an}) equals 1.
This proves (\ref{an}) with a ``$+$'' sign, and the same argument goes
through with ``$-$'' in place of ``$+$''.
\end{longlist}

(ii)
We use the It\^o representation in (\ref{decomp3}). 
Fix $t>0$ and take any constants $d_\pm>0$ satisfying~(\ref{Pilb}).
Let $\kappa_\pm>0$ be any constants and choose $K_\pm$ to satisfy
(\ref{Kpm}).
For the small jump processes, we have the bounds in (\ref{an}).
Note that these remain true if $\overline{\Pi}^+(0+)=0$ or $\overline
{\Pi}^-(0+)=0$.
For the big positive jumps, we have
\begin{eqnarray}
P\bigl(X_t^{(B,d_+,+)}=0\bigr) &\ge& P(\mbox{no }
\Delta X_s  \mbox{ exceeds } d_+\mbox{ up\ till time } t)
\nonumber
\\
\label{sldiff}
&=& e^{-t\overline{\Pi}^+(d_+)}
\\
&\ge & e^{-c_+}  \qquad \bigl(\mbox{by }(\ref{Pilb})\bigr).\nonumber
\end{eqnarray}
Equation (\ref{sldiff}) remains true with $c_+=0$ when $\overline{\Pi}^+(0+)=0$.
By (\ref{decomp3}), the probability on the left-hand
side of (\ref{Xlb2}) is, for any $L\ge0$,
\begin{eqnarray}
&& P\bigl(\sigma Z_t+ X_t^{(S,d_+,+)}+X_t^{(B,d_+,+)}+
X_t^{(S,d_-,-)}+X_t^{(B,d_-,-)}
\nonumber
\\
&& \qquad {}\le K_+d_+ -Ld_- -\kappa_+\sqrt{tV_+(d_+)} -\kappa_-\sqrt{tV_-(d_-)} \bigr)
\nonumber
\\
\label{ac4}
&&\quad\ge P \bigl(Z_t\le0, X_t^{(S,d_+,+)}\le K_+d_+-
\kappa_+\sqrt{tV_+(d_+)},\ X_t^{(B,d_+,+)}=0,
\\
&&\qquad  X_t^{(S,d_-,-)}\le K_-d_--\kappa_-\sqrt{tV_-(d_-)},\
X_t^{(B,d_-,-)}\le-(K_-+L)d_- \bigr)
\nonumber
\\
&&\quad \ge e^{-c_+} \Phi(-\kappa_+) \Phi(-\kappa_-) P \bigl(
X_t^{(B,d_-,-)}\le-(K_-+L)d_- \bigr)/8.\nonumber
\end{eqnarray}
In the last inequality,
we used (\ref{an}) (twice; once with ``$+$'' and once with ``$-$''),
(\ref{sldiff})
and the independence of the $Z_t$ and the $X_t^{(\cdot)}$ processes.
No jump in $X_t^{(B,d_-,-)}$ is larger than $-d_-$, so
we have the upper bound
$X_t^{(B,d_-,-)}\le-d_-N_t^-(d_-)$, where $N_t^-(d_-)$ is the number
of jumps of $X_t$ less than or equal in size to $-d_-$ which occur by
time $t$. $N_t^-(d_-)$ is distributed as Poisson with expectation
$t\overline{\Pi}^-(d_--)$, and $t\overline{\Pi}^-(d_--)\ge c_-$ by
(\ref{Pilb}).
(Note that this implies $\overline{\Pi}^-(0+)>0$.)
The Poisson distribution is stochastically monotone in the sense
that if $N(\mu_1)$ and $N(\mu_2)$
are Poisson rvs with means $\mu_1>\mu_2$, then
$P(N(\mu_1)\ge x) \ge P(N(\mu_2)\ge x)$ for all $x\ge0$.
So, letting $N(c_-)$
be a Poisson rv with expectation $c_-$, we have
\begin{equation}
\label{nc}
P\bigl(N_t^-(d_-)\ge K_-+L\bigr) \ge P\bigl(N(c_-)\ge
K_-+L\bigr).
\end{equation}
%
Then using
\begin{equation}
\label{ncn}
P \bigl( X_t^{(B,d_-,-)}\le-(K_-+L)d_- \bigr) \ge
P \bigl(N_t^-(d_-)\ge K_-+L \bigr)
\end{equation}
and (\ref{ac4}) we arrive at (\ref{Xlb2}).
When $\overline{\Pi}^+(0+)=0$, we can take all the ``$+$'' terms in
(\ref{ac4}) as 0 to get (\ref{Xlb2}) with all the ``$+$'' terms 0.

(iii)   Assume $0\le\overline{\Pi}^-(0+)<\infty=\overline{\Pi}^+(0+)$.
In this case, we do not define $d_-$ but still have $d_+=d_+(t)>0$ and
assume $t\overline{\Pi}^+(d_+)\le c_+$ as in (\ref{Pilb}). From
(\ref{decomp3}), write
\begin{equation}
\label{decomp5}
X_t=t\gamma-t\nu_+(d_+)+t\nu_-(0)
+X_t^{(S,d_+,+)} +X_t^{(B,d_+,+)} +
X_t^{(0,-)},
\end{equation}
where the negative jump components have been amalgamated into
\begin{eqnarray*}
&& X_t^{(0,-)}:= \sum_{0<s\le t} \Delta
X_s 1_{\{\Delta X_s\le0\}}, \qquad t>0,
\end{eqnarray*}
which is a compound Poisson process comprised of non-positive jumps.
This term and the term $t\nu_-(0)$ are absent from (\ref{decomp5})
when $\overline{\Pi}^-(0+)=0$.
Using (\ref{an}), (\ref{sldiff}) and (\ref{decomp5}), write
\begin{eqnarray*}
&& P \bigl(X_t\le t\gamma-t\nu_+(d_+)+t\nu_-(0+) +K_+d_+-
\kappa_+\sqrt{tV_+(d_+)} \bigr)
\nonumber
\\
&&\quad\ge P \bigl(Z_t\le0, X_t^{(S,d_+,+)}\le K_+d_+-
\kappa_+\sqrt{tV_+(d_+)}, X_t^{(B,d_+,+)}=0, X_t^{(0,-)}
\le0 \bigr)
\nonumber
\\
&&\quad\ge e^{-c_+}\Phi(-\kappa_+) P \bigl( X_t^{(0,-)}
\le0 \bigr)/4= e^{-c_+} \Phi(-\kappa_+)/4
\end{eqnarray*}
and this gives (\ref{Xlb3}).
\end{pf}

\section{Proof of Theorem~\texorpdfstring{\protect\ref{th3.3s}}{2.1}}\label{s2ac}
\textit{Part} (i).  Assume $\overline{\Pi}^-(0+)>0$ throughout
this part.

(\ref{3.6as}) $\Longrightarrow$ (\ref{3.7s}):  Assume (\ref{3.6as}).
$\overline{\Pi}^-(0+)>0$ implies $\overline{\Pi}^-(x)>0$ in a
neighbourhood of 0 so
we can assume $\overline{\Pi}^-(x)>0$ for all $0<x<1$.
Choose $1>x_k\downarrow0$ such that
\[
\frac{A(x_k)}{\sqrt{U(x_k)\overline{\Pi}^-(x_k)}}\to\infty
\]
as $k\to\infty$. This implies $\sigma^2=0$ by (\ref{4.1x})
(because $U(x)\ge\sigma^2$). 
It also means that $A(x_k)>0$ for all large $k$, and without loss of
generality we may assume it to be so for all $k$.
Let
\[
s_k:= \sqrt{\frac{U(x_k)}{\overline{\Pi}^-(x_k)A^2(x_k)}},
\]
then
\[
s_k\overline{\Pi}^-(x_k) = \frac{\sqrt{U(x_k)\overline{\Pi}^-(x_k)}}{A(x_k)}\to0
\]
and since $\overline{\Pi}^-(0+)>0$, also $s_k\to0$ as $k\to\infty$.
In addition, we have
\[
\frac{U(x_k)}{s_kA^2(x_k)} = \frac{\sqrt{U(x_k)\overline{\Pi}^-(x_k)}}{A(x_k)}\to0
\]
and
\[
\frac{s_kA(x_k)}{x_k}= \sqrt{\frac{U(x_k)}{x_k^2\overline{\Pi}^-(x_k)}}\ge1.
\]
Set
\[
t_k:= \sqrt{\frac{s_k}{\overline{\Pi}^-(x_k)}},
\]
so $t_k/s_k\to\infty$, but still
$t_k\overline{\Pi}^-(x_k)\to0$, as $k\to\infty$.
Then
\begin{equation}
\label{Ut0}
\frac{U(x_k)}{t_kA^2(x_k)} =\frac{s_k}{t_k}\frac{U(x_k)}{s_kA^2(x_k)} \to0,
\end{equation}
and
\begin{equation}
\label{Ati}
\frac{t_kA(x_k)}{x_k} =\frac{t_k}{s_k}\frac{s_kA(x_k)}{x_k} \to
\infty,
\end{equation}
as $k\to\infty$.

Recall (\ref{IP}) and use the It\^o decomposition in (\ref{decomp3})
with $\sigma^2=0$ and $h_+=h_-=h>0$ to write
\begin{equation}
\label{dec2} X_t=tA(h)+X_t^{(S,h)}+X_t^{(B,h,+)}-th
\overline{\Pi}^+(h)+ X_t^{(B,h,-)}+th\overline{\Pi}^-(h), \qquad
t>0.
\end{equation}
Here, $X_t^{(S,h)}=X_t^{(S,h,+)}+X_t^{(S,h,-)}$ is the compensated
small jump process, and $X_t^{(B,h,\pm)}$ are the positive and
negative big jump processes.
\begin{longlist}[\textit{Case} (a)]
\item[\textit{Case} (a):]  Suppose $\overline{\Pi}^+(0+)>0$.
Since each jump in $X_t^{(B,h,+)}$ is at least $h$,
we have the lower bound $X_t^{(B,h,+)}\geq hN_t^+(h)$,
where $N_t^+(h)$ is Poisson with expectation $t\overline{\Pi}^+(h)$
(and variance $t\overline{\Pi}^+(h)$). Using this and substituting in
(\ref{dec2})
with $t=t_k$ and $h=x_k$ we get
\begin{equation}
\label{dec4} X_{t_k} \ge t_kA(x_k)
+X_{t_k}^{(S,x_k)}+x_k \bigl(N^+_{t_k}(x_k)-t_k
\overline{\Pi }^+(x_k) \bigr) + X_{t_k}^{(B,x_k,-)}.
\end{equation}
Since $t_k\overline{\Pi}^-(x_k)\to0$,
we have $P(X_{t_k}^{(B,x_k,-)}=0)\to1$
as $k\to\infty$.
Also, for $\varepsilon\in(0,1)$,
\begin{eqnarray*}
P \bigl( X_{t_k}^{(S,x_k)}+x_k \bigl(N^+_{t_k}(x_k)-t_k
\overline{\Pi }^+(x_k) \bigr) \le-\varepsilon t_kA(x_k)
\bigr) &\le& \frac{t_kV(x_k)+t_kx_k^2\overline{\Pi}^+(x_k)}{\varepsilon^2t_k^2A^2(x_k)}
\nonumber
\\
&\le& \frac{U(x_k)}{\varepsilon^2t_kA^2(x_k)} \to0,
\end{eqnarray*}
as $k\to\infty$ by (\ref{Ut0}), so
\begin{equation}
\label{Xot}
P \biggl(\frac{X_{t_k}}{x_k} \ge(1-\varepsilon)\frac
{t_kA(x_k)}{x_k}
\biggr) \to1,
\end{equation}
and hence, by (\ref{Ati}), $X_{t_k}/x_k\stackrel{\mathrm
{P}}{\longrightarrow}\infty$ as $k\to\infty$.
Thus, (\ref{3.7s}) holds.

\item[\textit{Case} (b):]  Alternatively, if $\overline{\Pi}^+(0+)=0$, we can
omit the
term containing $N^+_{t_k}(x_k)-t_k\overline{\Pi}^+(x_k)$ in (\ref
{dec4}) and
in what follows it, and again obtain (\ref{Xot}), and hence (\ref
{3.7s}).\footnote{Observe that the assumption $\Pi(\mathbb
{R})=\infty$ was
not used in this part of the proof. The trivial case, $X_t=t\gamma$,
$\gamma>0$, when $A(x)\equiv\gamma$, is included if we interpret
(\ref{3.6as}) as holding then.}

(\ref{3.6as}) $\Longrightarrow$ (\ref{Xdels}):
Continuing the previous argument, $t_k\overline{\Pi}^-(x_k)\to0$ implies
\[
P \bigl( \bigl(\Delta X^-\bigr)_{t_k}^{(1)} >x_k
\bigr) = P \Bigl( \sup_{0<s\le t_k}\Delta X_s^-
>x_k \Bigr)= 1-e^{-t_k \overline{\Pi}^-(x_k)} \to0 \qquad \mbox{as } k\to\infty,
\]
so, using (\ref{Xot}), (\ref{Xdels}) also holds when (\ref{3.6as})
holds and $\overline{\Pi}^-(0+)>0$.

(\ref{Xdels}) $\Longrightarrow$ (\ref{3.7s}):
This is obvious when $\overline{\Pi}^-(0+)>0$.

(\ref{3.7s}) $\Longrightarrow$ (\ref{3.6as}):
Assume $\Pi(\mathbb{R})=\infty$ as well as $\overline{\Pi}^-(0+)>0$,
and that (\ref{3.7s}) holds. Suppose (\ref{3.6as}) fails, so we can choose
$1<a<\infty$, $x_0>0$, such that
\begin{equation}
\label{Ab}
A(x) \le a\sqrt{U(x)\overline{\Pi}^-(x)},
\end{equation}
for all $0<x\le x_0$.
We will obtain a contradiction.
Note that (\ref{3.7s}) implies $\sigma^2=0$, because $X_t/\sqrt
{t}\stackrel{\mathrm{D}}{\longrightarrow}N(0,\sigma^2)$, a
non-degenerate normal rv, when $\sigma^2>0$.
So we assume $\sigma^2=0$ in what follows.
We consider 3 cases.

\item[\textit{Case} (a):]  Assume in fact that $\overline{\Pi
}^-(0+)=\infty
=\overline{\Pi}^+(0+)$.
In this situation, we can introduce quantile versions for the $d_\pm$
in (\ref{Pilb}).
Define the non-decreasing function
\begin{equation}
\label{ddef}
d_+(t):=\inf\bigl\{x>0: \overline{\Pi}^+(x)\le t^{-1}
\bigr\},\qquad t>0,
\end{equation}
and set $d_+(0)=0$.
Since $\overline{\Pi}^+(0+)=\infty$, we have $0<d_+(t)<\infty$ for
all $t>0$,
$d_+(t)\downarrow0$ as $t\downarrow0$, and 
\begin{equation}
\label{dprop}
t\overline{\Pi}^+\bigl(d_+(t)\bigr) \le1\le t\overline{\Pi}^+
\bigl(d_+(t)-\bigr)\qquad \mbox{for  all } t>0.
\end{equation}
Analogously, define $d_-(0)=0$, and
\begin{equation}
\label{d-def} d_-(t):=\inf\bigl\{x>0: \overline{\Pi}^-(x)\le t^{-1}
\bigr\}, \qquad t>0,
\end{equation}
having, since $\overline{\Pi}^-(0+)=\infty$,
$0<d_-(t)<\infty$, $d_-(t)\downarrow0$ as $t\downarrow0$, and
\begin{equation}
\label{d-prop}
t\overline{\Pi}^-\bigl(d_-(t)\bigr) \le1 \le t\overline{\Pi}^-
\bigl(d_-(t)-\bigr).
\end{equation}
%

With $a$ as in (\ref{Ab}), set $\kappa_+=\kappa_-=\kappa=2a$, then
choose $K_\pm$ to satisfy (\ref{Kpm}).
Then (\ref{3.7s}) together with (\ref{Xlb2}) shows that we must have
\begin{equation}
\label{Ac}
0\le t_k \bigl(\gamma -\nu_+(d_+)+\nu_-(d_-) \bigr)
+K_+d_+-Ld_- -\kappa\sqrt{t_k\bigl(V_+(d_+)+ V_-(d_-)\bigr)},
\end{equation}
for all large $k$.
Here, $d_+$ and $d_-$ are any positive numbers
and we used the inequality
$\sqrt{a} +\sqrt{b} \ge\sqrt{a+b}$, $a,b>0$, in (\ref{Xlb2}).
Take $\lambda>0$ and set
\[
d_+= d_+(\lambda t_k)\quad \mbox{and}\quad d_-=d_-(t_k),
\]
where $d_+(\cdot)$ and $d_-(\cdot)$ are defined in (\ref{ddef}) and
(\ref{d-def}).
By (\ref{dprop}) and (\ref{d-prop}), we then have
\[
t\overline{\Pi}^+ \bigl(d_+(\lambda t)\bigr)\le\lambda^{-1}\quad\mbox{and} \quad t
\overline{\Pi}^-\bigl(d_-(t)-\bigr)\ge1,
\]
so (\ref{Pilb}) holds with $c_+=\lambda^{-1}$ and $c_-=1$. With $t_k$
as the sequence in (\ref{3.7s}), let $d=d(t_k):=\max(d_+(\lambda
t_k),d_-(t_k))$.

Equation (\ref{Ac}) implies
\begin{eqnarray}
&& t_k \biggl(\gamma-\int_{(d_+,1]}y\Pi(\mathrm{d}y)+\int_{(d_-,1]}y\Pi ^{(-)}(\mathrm{d}y)
\biggr)
\nonumber
\\[-8pt]
\label{Ac0}
\\[-8pt]
\nonumber
&&\quad\ge -K_+d_+ + Ld_- +\kappa\sqrt{t_k\bigl(V_+(d_+)+ V_-(d_-)
\bigr)}.
\end{eqnarray}
Adding the quantity
\[
t_k \biggl(\int_{(d_+,d]}y\Pi(\mathrm{d}y)-\int
_{(d_-,d]}y\Pi ^{(-)}(\mathrm{d} y) \biggr)
\]
to both sides of (\ref{Ac0}) gives $t_k\nu(d)$ on the left, and a
quantity no smaller than
\[
t_kd_+ \bigl(\overline{\Pi}^+(d_+)-\overline{\Pi}^+(d) \bigr)
-t_kd \bigl(\overline{\Pi}^-(d_-)-\overline{\Pi}^-(d) \bigr) -K_+d_+
+ Ld_- +\kappa\sqrt{t_k\bigl(V_+(d_+)+ V_-(d_-)\bigr)}
\]
on the right. Further adding $t_kd(\overline{\Pi}^+(d)-\overline{\Pi}^-(d))$
to both sides gives $t_kA(d)$ on the left
(see~(\ref{IP})),
and then after some cancellation we arrive at
\begin{equation}
\label{Ac1} t_kA(d)\ge t_kd_+\overline{
\Pi}^+(d_+)-t_kd\overline{\Pi}^-(d_-) -K_+d_+ + Ld_- + \kappa
\sqrt{t_k\bigl(V_+(d_+)+ V_-(d_-)\bigr)}.
\end{equation}

At this stage, it is helpful to assume that $\overline{\Pi}^+(x)$ is a
continuous function on $(0,\infty)$. It then follows from (\ref
{dprop}) that
$t_k\overline{\Pi}^+(d_+(\lambda t_k))=1/\lambda$, while
$t_k\overline{\Pi}^-(d_-(t_k))\le1$ by (\ref{d-prop}).
Also, $d\le d_++d_-$.
Thus, we deduce
\begin{equation}
\label{Ad}
t_kA(d)\ge (1/\lambda-K_+- 1 )d_+ + (L- 1 )d_- +
\kappa\sqrt{t_k\bigl(V_+(d_+)+ V_-(d_-)\bigr)}.
\end{equation}
Next, write
\begin{eqnarray}
V_+(d_+)+ V_-(d_-) &=& V_+(d)-\int_{(d_+,d]}y^2
\Pi(\mathrm{d}y) +V_-(d)-\int_{(d_-,d]}y^2
\Pi^{(-)}(\mathrm{d}y)
\nonumber
\\
\label{Vest}
&\ge& V(d)-d^2 \bigl(\overline{\Pi}^+(d_+)-\overline{\Pi}^+(d)
\bigr) -d^2 \bigl(\overline{\Pi}^-(d_-)-\overline{\Pi}^-(d) \bigr)
\\
&=& U(d)-d^2 \bigl(\overline{\Pi}^+(d_+)+\overline{\Pi}^-(d_-)
\bigr).\nonumber
\end{eqnarray}
So
\[
t_k \bigl( V_+(d_+)+ V_-(d_-) \bigr) \ge t_kU(d)-d^2
(1/\lambda+1 )
\]
giving
\[
\sqrt{t_k\bigl(V_+(d_+)+ V_-(d_-)\bigr)} \ge \sqrt{t_kU(d)}-d
(1/\sqrt{\lambda}+1 ).
\]
Substituting into (\ref{Ad}), we obtain
\begin{eqnarray*}
t_kA(d) &\ge& \kappa\sqrt{t_kU(d)} + (1/\lambda-K_+- 1 -
\kappa/\sqrt{\lambda}-\kappa )d_+
\nonumber
\\
&& {}+ (L-1-\kappa/\sqrt{\lambda}-\kappa )d_-.
\end{eqnarray*}
Choose $\lambda$ small enough for the
first expression in brackets on the right-hand side to be positive.
Then choose $L$ large enough for the
second expression in brackets on the right-hand side to be positive.
This gives
\begin{equation}
\label{Ag}
t_kA(d)\ge\kappa\sqrt{t_kU(d)},
\end{equation}
for all large $k$, where $d=d(t_k)= \max(d_+(\lambda
t_k),d_-(t_k))\downarrow0$ as $k\to\infty$. Inequality
(\ref{Ag}) implies
\begin{equation}
\label{bn}
\frac{A(d)}{\sqrt{U(d)\overline{\Pi}^-(d)}} \ge \frac{\kappa}{\sqrt{t_k\overline{\Pi}^-(d)}} \ge \frac{\kappa}{\sqrt{t_k\overline{\Pi}^-(d_-)}}
\ge\kappa,
\end{equation}
giving a contradiction with (\ref{Ab}), since $\kappa=2a$.

This proves (\ref{3.6as}) from (\ref{3.7s}) in case $\overline{\Pi}
^+(0+)=\overline{\Pi}^-(0+)=\infty$ and $\overline{\Pi}^+(x)$ is
continuous for $x>0$.
To complete the proof of part (i), case (a), of the theorem we remove
the assumption of continuity made in deriving (\ref{Ad}). This can be
done using the following lemma.
\end{longlist}

\begin{lemma}\label{lemexvaguenicepis} Let
$\Pi$ be any L\'evy measure with $\overline{\Pi}^+(0+)=\infty$.
Then there
exists a sequence of L\'evy measures $\Pi_n$,
absolutely continuous with respect to Lebesgue measure,
and having strictly positive $C^\infty$-densities on $\mathbb
{R}_*:=\mathbb{R}
\setminus\{0\}$,
satisfying $\overline{\Pi}_n^+(0+)=\infty$ and $\Pi_n\stackrel
{v}{\longrightarrow}\Pi$ as $n\to
\infty$.
\end{lemma}

\begin{pf}
($\stackrel{v}{\longrightarrow}$ refers to vague convergence in
$\overline\mathbb{R}_*$; see,
for example, Chapter~15 in Kallenberg \cite{Kall}.)
We extend $\Pi$ to a Borel measure on $\mathbb{R}$ by setting $\Pi(\{
0\}
):=0$. Assume $\Pi\neq0$, so~$C:=\int x^2 \Pi(\mathrm{d}x)/(1+x^2)\in(0,\infty)$. Observe that
$P(\mathrm{d}x):=x^2 \Pi(\mathrm{d}x)/C(1+x^2)$ defines
a Borel probability measure on $\mathbb{R}$. For all $n\in\mathbb
{N}$, the
convolved probability measure $P_n:=P\star N(0,1/n)$ admits a strictly positive
$C^\infty$-Lebesgue density, when $N(0,1/n)$ is a normal rv with
expectation 0 and variance $1/n$. Set $\Pi_n(\mathrm{d}x):= C(1+x^2)
P_n(\mathrm{d}x)/x^2$, $n=1,2,\ldots.$
It is easily verified that $(\Pi_n)_{n\in\mathbb{N}}$ is a sequence
of L\'
evy measures with the desired properties.
\end{pf}

Now to complete the proof of part (i), assume (\ref{3.7s}), that is,
$P(X_{t_k}\ge0)\to1$ for some $t_k\downarrow0$, for a general $X$
with L\'
evy measure $\Pi$.
Using Lemma~\ref{lemexvaguenicepis}, construct a sequence of
approximating L\'evy measures $\Pi_n$, converging vaguely to $\Pi$,
such that their positive tails $\overline{\Pi}_n^+$
are continuous on $(0,\infty)$ with $\overline{\Pi}_n^+(0+)=\infty$.
On the negative side, let $\overline{\Pi}_n^-(x)=\overline{\Pi
}^-(x)$, $x>0$.
Let $(X_t(n))_{t\ge0}$ be L\'evy processes with measures $\Pi_n$ and
other characteristics the same as for $X$. Define $\nu_n$, $A_n$,
$V_n$, $U_n$,
$\nu_{n,\pm}$, $V_{n,\pm}$, as in (\ref{Adef}), (\ref{Vdef}),
(\ref{nupmdef}) and (\ref{Vpmdef}), but with $\Pi_n$ replacing $\Pi
$. The subscript $n$ functions converge to the original functions at
points of continuity of the latter.
$X_t(n)$ has characteristic exponent given by (\ref{ce}) with $\Pi_n$
replacing $\Pi$, so as
$n\to\infty$ we have $X_t(n)\stackrel{\mathrm{D}}{\longrightarrow
}X_t$ for each $t>0$.
Under assumption (\ref{3.7s}), $\lim_{n\to\infty}P(X_{t_k}(n)>0)
=P(X_{t_k}>0)>1-\delta$ for arbitrary $\delta\in(0,1/2)$ and $k$
large enough. Thus, $P(X_{t_k}(n)>0)>1-2\delta$ for $n\ge n_0(k)$ and
$k\ge k_0$. So (\ref{Ac}) holds
for the subscript $n$ quantities
with probability $1-2\delta$. But (\ref{Ac}) is deterministic so it
holds in fact for the subscript $n$ quantities (with probability 1)
whenever $n\ge n_0(k)$ and $k\ge k_0$.
The proof using continuity of $\Pi_n$ then shows that
(\ref{bn}) holds with $A$, $U$ and $\overline{\Pi}^-$ replaced by $A_n$,
$U_n$ and $\overline{\Pi}_n^-$.
Then letting $n\to\infty$ shows that (\ref{bn}) itself holds as
stated for $k\ge k_0$.
Again we get a contradiction, and thus complete the proof that (\ref
{3.7s}) implies (\ref{3.6as}) for case (a).
\begin{longlist}[\textit{Case} (b).]
\item[\textit{Case} (b):]  Assume that $0\le\overline{\Pi
}^+(0+)<\infty
=\overline{\Pi}^-(0+)$.
As in the proof for case (a), we take $\kappa_+=\kappa_-=\kappa>2a$,
$K_\pm$ to satisfy (\ref{Kpm}), define $d_-(t)>0$ by (\ref{d-def}),
and write $d_-=d_-(t)$ for $t>0$.
But for $d_+$ we set $d_+(t)\equiv d_-(t)>0$.
We take $c_+=1$ in (\ref{Pilb}) as we may since $t\overline{\Pi}
^+(d_+)=t\overline{\Pi}^+(d_-)\le t\overline{\Pi}^+(0+)\to0$ as
$t\downarrow0$.
With this set-up, (\ref{Ac}) is true (with $d_-$ replacing $d_+$) and
we can follow the proof of case (a) through to get (\ref{Ac1}) with
$d$ and $d_+$ replaced by $d_-$; thus,
\begin{equation}
\label{Ac10} t_kA(d_-)\ge t_kd_- \bigl(\overline{
\Pi}^+(d_-)-\overline{\Pi}^-(d_-) \bigr) -K_+d_- + Ld_- + \kappa
\sqrt{t_k\bigl(V_+(d_-)+V_-(d_-)\bigr)}.
\end{equation}
Estimating $V_\pm$ along the lines of (\ref{Vest}) we find the
right-hand side of (\ref{Ac10}) is not smaller than
\begin{eqnarray*}
&&\kappa\sqrt{t_kU(d_-)} + (L-K_+-1-\kappa )d_-.
\end{eqnarray*}
Choose $L$ large enough in this to get (\ref{Ag}), and hence
(\ref{bn}) with $d_-$ in place of $d$, hence (\ref{3.6as}) again.

\item[\textit{Case} (c):]  Assume that $0<\overline{\Pi
}^-(0+)<\infty=\overline{\Pi}
^+(0+)$. Define $d_+(t)$ by (\ref{ddef}), so we have (\ref{dprop}).
Then (\ref{Xlb3}) with $c_+=1$ and $\kappa_+=\kappa=2a$, together
with (\ref{3.7s}), shows that we must have
\begin{equation}
\label{Ac6}
0\le t_k \bigl(\gamma-\nu_+(d_+)+\nu_-(0)
\bigr)+K_+d_+ -\kappa\sqrt{t_kV_+(d_+)},
\end{equation}
for all large $k$.
Here again we write $d_+= d_+(\lambda t_k)$
for $\lambda>0$.
Inequality (\ref{Ac6}) implies
\begin{equation}
\label{Ac7}
t_k \biggl(\gamma-\int_{(d_+,1]}y\Pi(\mathrm{d}y)+\int_{(0,1]}y\Pi ^{(-)}(\mathrm{d}y)
\biggr) \ge -K_+d_+ +\kappa\sqrt{t_kV_+(d_+)}.
\end{equation}
Subtracting the quantity
\[
t_k\int_{(0,d_+]}y\Pi^{(-)}(\mathrm{d}y)\le
t_kd_+ \bigl(\overline {\Pi} ^-(0+)-\overline{\Pi}^-(d_+) \bigr)
\]
from both sides of (\ref{Ac7}) gives $t_k\nu(d_+)$ on the left, and a
quantity no smaller than
\[
-t_kd_+ \bigl(\overline{\Pi}^-(0+)-\overline{\Pi}^-(d_+) \bigr)
-K_+d_+ +\kappa\sqrt{t_kV_+(d_+)}
\]
on the right. Further adding $t_kd_+(\overline{\Pi}^+(d_+)-\overline
{\Pi}^-(d_+))$
to both sides gives (see (\ref{IP}))
\begin{eqnarray}
\label{Ac8} t_kA(d_+)\ge t_kd_+\overline{
\Pi}^+(d_+)-t_kd_+\overline{\Pi}^-(0+)-K_+d_+ + \kappa
\sqrt{t_kV_+(d_+)}.
\end{eqnarray}

At this stage, as before, assume $\overline{\Pi}^+(x)$ is continuous.
It then
follows from (\ref{dprop}) that
$t_k\overline{\Pi}^+(d_+(\lambda t_k))=1/\lambda$, while
$t_k\overline{\Pi}^-(0+)\le1$ for large $k$.
Thus, from (\ref{Ac8}) we deduce
\begin{equation}
\label{Ad1}
t_kA(d_+)\ge (1/\lambda-K_+- 1 )d_+ +\kappa
\sqrt{t_kV_+(d_+)},
\end{equation}
for large enough $k$. Further,
\begin{eqnarray*}
t_kV_+(d_+)&=& t_k \bigl(U(d_+)-V_-(d_+)-d_+^2
\overline{\Pi}^+(d_+)-d_+^2\overline {\Pi}^-(d_+) \bigr)
\cr
&\ge&
t_kU(d_+)-t_kd_+^2 \bigl(2\overline{
\Pi}^-(0+)+\overline{\Pi }^+(d_+) \bigr)
\cr
&\ge& t_kU(d_+)-4d_+^2,
\end{eqnarray*}
using that $V_-(d_+)\le d_+^2\overline{\Pi}^-(0+)$. So
\[
\sqrt{t_kV_+(d_+)}\ge\sqrt{t_kU(d_+)}-2d_+.
\]
Substituting into (\ref{Ad1}), we get
\begin{eqnarray*}
t_kA(d_+) \ge \kappa\sqrt{t_kU(d_+)}+ (1/\lambda-K_+-1-2
\kappa )d_+,
\end{eqnarray*}
for large $k$.
We can choose $\lambda$ small enough for the expression in brackets on
the right-hand side to be positive.
This gives
\begin{equation}
\label{Ac9}
t_kA(d_+)\ge\kappa\sqrt{t_kU(d_+)},
\end{equation}
which, since $\overline{\Pi}^-(0+)>0$ is assumed, implies
\begin{eqnarray*}
&& \frac{A(d_+)}{\sqrt{U(d_+)\overline{\Pi}^-(d_+)}} \ge \frac{\kappa}{\sqrt{t_k\overline{\Pi}^-(0+)}} \to\infty \qquad \mbox{as } k\to\infty,
\end{eqnarray*}
a contradiction with (\ref{Ab}). We can remove the continuity
assumption as before.
So (\ref{3.6as}) is proved when $0<\overline{\Pi}^-(0+)<\infty
=\overline{\Pi}^+(0+)$.

\item[\textit{Part} (ii).]
Now we will deal with the case when
$\overline{\Pi}^-(0+)=0$ but $\overline{\Pi}^+(0+)=\infty$.
Again assume $\overline{\Pi}^+(x)$ is continuous.
The working in case (c) is still valid from (\ref{Ac6}) to (\ref{Ac9}).
The negative jump process is now absent from $X_t$ and
(\ref{Ac6}) gives, with $d_+=d_+(\lambda t)$,
\begin{eqnarray}
0&\le& t_k \bigl(\gamma-\nu_+(d_+) \bigr)+K_+d_+\nonumber
\\
\label{Ac2}
&=& t_k \biggl(\gamma-\int_{d_+}^1
\overline{\Pi}^+(y)\,\mathrm{d}y -d_+\overline{\Pi}^+(d_+)+\overline{\Pi}^+(1)
\biggr)+K_+d_+
\\
&\le& t_k \biggl(\gamma-\int_{d_+}^1
\overline{\Pi}^+(y)\,\mathrm{d}y +\overline{\Pi}^+(1) \biggr)-(1/\lambda-K_+)d_+,\nonumber
\end{eqnarray}
using $t_k\overline{\Pi}^+(d_+(\lambda t_k))= 1/\lambda$ in the last
inequality
(since $\overline{\Pi}^+(x)$ is continuous).
Then choosing $\lambda<1/K_+$ we get
\begin{eqnarray*}
&&\int_{d_+}^1\overline{\Pi}^+(y)\,\mathrm{d}y \le
\gamma+\overline {\Pi}^+(1).
\end{eqnarray*}
Letting $k\to\infty$ (so $t_k\downarrow0$ and $d_+=d_+(\lambda
t_k)\downarrow
0$) shows that $\int_0^1\overline{\Pi}^+(y)\,\mathrm{d}y<\infty$.
Since $X_t$ has no negative jumps, from this we deduce that $X$ is of
bounded variation with drift $\mathrm{d}_X=A(0+)$
(see Doney and Maller \cite{dm1}, Theorem~2.1 and Remark~1), which is
non-negative by (\ref{Ac9}). Thus, $X$ is a subordinator with
non-negative drift.
It follows that
\[
A(x)= \gamma+\overline{\Pi}^+(1)-\int_x^1
\overline{\Pi }^+(y)\,\mathrm{d}y = \mathrm{d}_X+\int
_0^x \overline{\Pi}^+(y)\,\mathrm{d}y
\]
is non-negative for all $x>0$.
This is proved assuming continuity of $\Pi_n$ but that assumption can
be removed as before.
Then $A(x)\ge0$ together with $\sigma^2=0$ implies $\lim_{t\downarrow
0}P(X_t>0)\to1$ as $t\downarrow0$ by Theorem~\ref{th3.3}, hence
(\ref{3.7s}).

\item[\textit{Part} (iii).]
Finally, suppose $\overline{\Pi}^-(0+)>0$ and
(\ref{3.6as+}) holds. Then we can choose $x_k\downarrow0$
such that
\[
\frac{A(x_k)}{\sqrt{U(x_k)\overline{\Pi}^-(x_k)}}\to\infty \quad \mbox{and}\quad A(x_k) \to\infty,
\]
as $k\to\infty$.
Following exactly the proof of part (i), we get (\ref{Xot}), and this
implies $X_{t_k}/t_k\stackrel{\mathrm{P}}{\longrightarrow}\infty$,
since $A(x_k)\to\infty$.

Conversely, suppose
$X_{t_k}/t_k\stackrel{\mathrm{P}}{\longrightarrow}\infty$ as $k\to
\infty$ for a non-stochastic
sequence $t_k\downarrow0$. Then $\lim_{k\to\infty}P(X_{t_k}^M>0)=1$ for
every $M>0$, where $X_t^M$ is L\'evy with triplet $(\gamma-M,\sigma
^2,\Pi)$.
Consequently, (\ref{3.6as}) holds with $A$, $U$, $\Pi$ replaced by
$A^M(\cdot)=A(\cdot)-M$, $U^M=U$, $\Pi^M=\Pi$,
and this modified version implies (\ref{3.6as+}).
This completes part (iii),
and the proof of the theorem.
\end{longlist}

\section*{Acknowledgements}
I am grateful to a
referee for a close reading of the paper and helpful suggestions, and
to Boris Buchmann for supplying Lemma~\ref{lemexvaguenicepis}.

Research partially supported by ARC Grant DP1092502.






\printhistory
\end{document}